\documentclass[12pt,a4paper]{article}
\usepackage{amssymb}
\usepackage{amsmath}
\usepackage{indentfirst}

\newtheorem{theorem}{Theorem}

\newtheorem{corollary}[theorem]{Corollary}

\newtheorem{definition}{Definition}

\newtheorem{proposition}[theorem]{Proposition}
\newtheorem{remark}{Remark}

\newenvironment{proof}[1][]{\textbf{Proof.} }{}

\begin{document}

\title{Lagrange--Fedosov Nonholonomic Manifolds}

\author{Fernando Etayo\thanks{Departamento de Matem\'{a}ticas, Estad\'{\i}stica
y Computaci\'{o}n, Facultad de Ciencias, Universidad de Cantabria,
Avda. de los Castros s/n, 39071 Santander SPAIN;  e-mail: etayof@unican.es}\ ,
Rafael Santamar\'{\i}a
\thanks{Departamento de Matem\'{a}ticas, Escuela de Ingenier\'{\i}as
Industrial e Inform\'{a}tica, Universidad de Le\'{o}n, 24071 Le\'{o}n; e-mail: demrss@unileon.es}\ \  and Sergiu
I. Vacaru\thanks{Instituto de Matem\'{a}ticas y F\'{i}sica Fundamental, Consejo Superior de Investigaciones
Cient\'{i}ficas, Calle Serrano 123, 28006 Madrid SPAIN; e--mail: vacaru@imaff.cfmac.csic.es,
sergiu\underline{\hspace{1mm}}vacaru@yahoo.com}
 \\
Departamento de Matem\'{a}ticas, Estad\'{\i}stica y Computaci\'{o}n,\\
Facultad de Ciencias, Universidad de Cantabria,\\
Avda. de los Castros, s/n, 39071 Santander, Spain}
\date{}
 \maketitle

\begin{abstract}

We outline an unified approach to geometrization of Lagrange mechanics,
Finsler geometry and geometric methods
of constructing exact solutions with generic off--diagonal terms and nonholonomic
variables in gravity theories.
Such geometries with induced almost
symplectic structure are modelled on nonholonomic manifolds provided with nonintegrable
distributions defining nonlinear connections. We introduce the concept of
Lagrange--Fedosov spaces and Fedosov nonholonomic manifolds provided with almost
symplectic connection adap\-ted to the nonlinear connection structure.
We investigate the main properties of generalized Fedosov nonholonomic manifolds and
analyze exact solutions defining almost symplectic Einstein spaces.

\vskip0.3cm AMS Subject Classification:

 51P05, 53D15,  53B40,
53C07, 53C55,  70G45, 83C15

\vskip0.3cm PACS Classification:\  02.40.Yy, 45.20Jj, 04.20.Jb, 04.90.+3

\vskip0.3cm \textbf{Keywords:}\ Fedosov, Lagrange and
Finsler geometry, nonlinear connection, almost symplectic and almost
product structures, nonholonomic manifolds, exact solutions in gravity.
\end{abstract}

\newpage

\section{ Introduction}

The geometry of Fedosov manifolds is a natural generalization of
K\"{a}hler geometry defining a procedure of canonical deformation
quantization$^{1-5}$. By definition, a Fedosov manifold is given by a triple $%
(M,\theta ,\Gamma )$ where $M$ is a $\mathcal{C}^{\infty
}$--manifold enabled with symplectic structure $\theta$ (a
non--degenerated closed exterior 2--form) and a symplectic
connection structure $\Gamma $ (i.e. a torsionless linear
connection parallelizing the symplectic form). If a Lagrange
fundamental function $L:(x,y)\in TM\rightarrow \mathbb{R}$ is
defined
on $M$ \footnote{For simplicity, in this work we shall consider only regular Lagrangians; $%
(x,y)$ denote a set of local coordinates on the tangent bundle $TM$ with $%
x\in M$.}, there is a natural almost complex structure adapted to
the canonical nonlinear connection (in brief, N--connection)
induced by $L(x,y)$$^{6-8}$. Nonlinear connections can be also
naturally related to generic off--diagonal metrics and
nonholonomic moving frames in (super) gravity and string
theories$^{9-13}$. So, if we want to apply the methods of
symplectic geometry (and possible generalizations for Poisson
manifolds$^{14}$) to various type of Lagrange--Hamilton, and
related Finsler--Cartan spaces, we have to consider spaces enabled
with N--connection structure.

In this work, we study the geometry of almost symplectic
connections (in general, they are not torsion--free but can be
symmetrized) which are distinguished by a N--connection structure
and preserve an almost symplectic form, for instance, induced by a
regular Lagrangian or off--diagonal metric structure. This is
related to almost symplectic manifolds (see, for
instance$^{15-19}$) but, in our case, the manifolds are
nonholonomic ones.

We shall define and analyze the curvature tensor for such almost
symplectic connections and related Einstein equations with
nonholonomic variables. For nonholonomic manifolds, i.e.,
manifolds with nonintegrable distributions (in our case, with a
such distribution defined by a N--connection), this is not a
trivial task. The problem together with a proposal when the
Riemann tensor is interpreted as a modification of the Spencer
cohomology and related to solutions of partial differential
equations, as well to superspaces, are analyzed in$^{20,21}$.

The geometry of nonholonomic manifolds has a long time historical
perspective: For instance, in the review$^{22}$ it is stated that
it is probably impossible to construct an analog of the Riemannian
tensor for the general nonholonomic manifold. In two more recent
reviews$^{23,24}$, it is emphasized that in the past there were
proposed well defined Riemannian tensors for a number of spaces
provided with nonholonomic distributions, like Finsler and
Lagrange spaces and various type of theirs higher order
generalizations, i.e., for nonholonomic manifolds possessing
corresponding N--connection structures. As some examples of former
such investigations, we cite the works$^{25-29}$.

Essentially, the Fedosov type nonholonomic geometry to be elaborated in this work is based on the notion of
N--connection and considers a Whitney-like splitting of the tangent bundle to a manifold into  horizontal and
vertical subspaces (see discussion and a bibliography for recent developments and applications in$^{30-32}$).
Here we emphasize that the geometrical aspects of the N--connection formalism has been studied since the first
papers of E. Cartan$^{33}$ and A. Kawaguchi$^{34-36}$ (who used it in component form for Finsler geometry), then
one should be mentioned the so called Ehresmann connection$^{37}$ and the work of W. Barthel$^{38}$ where the
global definition of N--connection was given. The monographs$^{6-8}$ consider the N--connection formalism
elaborated and applied to the geometry of generalized Finsler--Lagrange and Cartan--Hamilton spaces, see also
the approaches$^{39-42}$.

The works related to nonholonomic geometry and N--connections have
appeared many times in a rather dispersive way when different
school of authors from geometry, mechanics and physics have worked
many times not having relation with another. We outline some
recent results with explicit applications in modern mathematical
physics and particle and string theories: N--connection structures
were modelled on Clifford and spinor bundles$^{43,44}$, on superbundles and in some directions
 of (super) string theory$^{45,46}$, as well in noncommutative geometry and gravity$^{47}$. The idea to apply the N--connections formalism as a new
geometric method of constructing exact solutions in gravity
theories was suggested in$^{9,10}$ and developed in a number of
works, see for instance$^{11-13}$).

We begin in Section II with an introduction into the N--connection
geometry for arbitrary manifolds with tangent bundles admitting
splitting into conventional horizontal and vertical subspaces. We
illustrate how regular Lagrangians induce natural semispray,
N--connection, metric and almost complex structures on tangent
bundles and discuss the relation between Lagrange and Finsler
geometry and theirs generalizations. Then we prove that
N--connection structures and corresponding almost complex
geometries may be modelled by generic off--diagonal metrics and
nonholonomic frames in gravity theories.

Section III is devoted to the theory of linear connecti\-ons on
N--anholono\-mic manifolds (i.e., on manifolds with nonholonomic
structure defined by N--connections). We demonstrate how the
linear connections may be adapted to the N--connection splitting
of the manifolds and analyze the conditions when such
distinguished connections may be naturally related to almost
complex structures. This has great {\em philosophical} interest,
because several authors have  defined different notions of general
connections, looking for associated  parallel transport and
covariant differential operator   satisfying, if possible,  the
properties of those of a linear connection (e.g., Ehresmann
connections on  bundles, non-homogeneous connections of
Grifone$^{48}$, quasi-- and pseudo--connections (see the
survey$^{49}$, etc.), but it always implies to lose properties or
to demand more assumptions than in the case of a N--connection
\footnote{For example, non-homogeneous connections of Grifone
define a covariant derivative $D_{X}Y$, which in general does not
define a vector field on the manifold (see p. 302 of the
reference$^{48}$) and which does not satisfy $D_{X}(fY)=fD_{X}Y +
(Xf)Y$ (see p. 305 of the same reference). In the case of
nonlinear connections used in the book of Yano and Ishihara$^{50}$
p. 209 it is assumed that horizontal distributions are invariant
under dilatations  (see also$^{51}$), etc.}. In the present paper
we shall see that one can define a canonical linear connection
adapted to a given N--connection. This shall avoid us extra
constructions and additional restrictions.

In Section IV, we define the Fedosov N--anholonomic and
Lagrange--Fedos\-ov ma\-nifolds as certain generalizations of the
Fedosov spaces to non\-holonomic configurations. We construct in
explicit form the curvature tensor of such spaces and define the
Einstein equations for N--adapted linear connection and metric
structures.

In Section V, we analyze the main conditions when vacuum
gravitational configurations with N--anholonomic structures can be
defined as exact solutions of the Einstein equations. We prove
that for a very general five dimensional ansatz for metric
coefficients depending on two, three and four variables the system
of field equations is completely integrable. We illustrate that
the method can be reduced to the case of four dimensional spaces
which gives us the possibility to generate conformal almost
complex gravitational metrics.

  We shall use both physical and mathematical languages and both coordinate and
intrinsic  notations, when possible.

\section{Nonlinear Connections and Fedosov Spa\-ces}

In this section, we recall some results on nonlinear connections and almost symplectic structures, which in
certain particular cases, are induced by regular Lagrangians, Finsler fundamental functions or by generic
off--diagonal metrics in gravity theories. From now on, all the manifolds
(in general, nonholonomic ones) \footnote{%
In literature, it is also used the equivalent term: anholonomic.} and geometric objects are supposed to be
$\mathcal{C}^{\infty }.$

\subsection{Nonlinear connection geometry}

Let $\mathbf{V}$ be a $(n+m)$--dimensional manifold. It is supposed that in
any point $u\in \mathbf{V}$ there is a local splitting \ $\mathbf{V}%
_{u}=M_{u}\oplus V_{u},$ where $M$ is a $n-$dimensional subspace and $V$ is a $m$--dimensional subspace
\footnote{One has this local decomposition when $\mathbf{V}\to M $ is a surjective submersion. A particular case
is that of a fibre bundle, but we can obtain the results in the general case. }. We shall split the local
coordinates (in general, abstract ones both for holonomic and nonholonomic variables) in the
form $u=(x,y),$ or $u^{\alpha }=\left( x^{i},y^{a}\right) ,$ where $%
i,j,k, \ldots =1,2, \ldots ,n$ and $a,b,c, \ldots =n+1,n+2, \ldots ,n+m.$ We denote by
$\pi ^{\top }:T\mathbf{V}\rightarrow TM$ the
differential of a map $\pi :V^{n+m}\rightarrow V^{n}$ defined by fiber preserving morphisms of the tangent
bundles $T\mathbf{V}$ and $TM.$ The kernel of $\pi ^{\top }$is just
the vertical subspace $v\mathbf{V}$ with a related inclusion mapping $i:v%
\mathbf{V}\rightarrow T\mathbf{V}.$

\begin{definition}
A nonlinear connection (N--connection) $\mathbf{N}$ on a manifold $\mathbf{V} $ is defined by the splitting on
the left of an exact sequence
\begin{equation*}
0\rightarrow v\mathbf{V}\overset{i}{\rightarrow} T\mathbf{V}\rightarrow T%
\mathbf{V}/v\mathbf{V}\rightarrow 0,
\end{equation*}%
i. e. by a morphism of submanifolds $\mathbf{N:\ \ }T\mathbf{V}\rightarrow v%
\mathbf{V}$ such that $\mathbf{N\circ i}$ is the unity in $v\mathbf{V}.$
\end{definition}

In an equivalent form, we can say that a N--connection is defined by a splitting
to subspaces with a Whitney sum
of conventional horizontal (h)
subspace, $\left( h\mathbf{V}\right) ,$ and vertical (v) subspace, $\left( v%
\mathbf{V}\right) ,$
\begin{equation}
T\mathbf{V}=h\mathbf{V}\oplus v\mathbf{V}  \label{whitney}
\end{equation}%
where $h\mathbf{V}$ is isomorphic to $M.$ Moreover, one can say that a N--connection is defined by a tensor
field of type (1,1) $\mathbf{P}=\mathbf{H}-\mathbf{N}$, where $\mathbf{H}$ (resp. $\mathbf{N}$) denotes the
projection over the horizontal (resp. vertical) subspace. Observe that $\mathbf{P}\circ \mathbf{P}=\mathbf{I}$,
i.e., $\mathbf{P}$ is an almost product structure, horizontal (resp. vertical) subspace being the eigenspace
associated to the eigenvalue +1 (resp. -1).

Locally, a N--connection is defined by its coefficients $N_{i}^{a}(u),$%
\begin{equation*}
\mathbf{N}=N_{i}^{a}(u)dx^{i}\otimes \frac{\partial }{\partial y^{a}}.
\end{equation*}%
The well known class of linear connections consists on a particular subclass with the coefficients being linear
on $y^{a},$ i.e., $N_{i}^{a}(u)=\Gamma _{bj}^{a}(x)y^{b}.$

Any N--connection $\mathbf{N}=N_{i}^{a}(u)$ may be characterized by an associated frame (vielbein) structure
$\mathbf{e}_{\nu }=(e_{i},e_{a}),$ where
\begin{equation}
e_{i}=\frac{\partial }{\partial x^{i}}-N_{i}^{a}(u)\frac{\partial }{\partial y^{a}}\mbox{ and
}e_{a}=\frac{\partial }{\partial y^{a}},  \label{dder}
\end{equation}%
and the dual frame (coframe) structure $\mathbf{\vartheta }^{\mu }=(\vartheta ^{i},\vartheta ^{a}),$ where
\begin{equation}
\vartheta ^{i}=dx^{i}\mbox{ and }\vartheta ^{a}=dy^{a}+N_{i}^{a}(u)dx^{i}. \label{ddif}
\end{equation}%
These vielbeins are called N--adapted frames. In order to preserve
a relation with the previous denotations$^{8-11,43,44,46}$, we
note that $\mathbf{e}_{\nu }=(e_{i},e_{a})$ and $\mathbf{\vartheta
}^{\mu }=(\vartheta ^{i},\vartheta ^{a})$ are, respectively, the
former $\delta _{\nu }=\delta /\partial u^{\nu }=(\delta
_{i},\partial _{a})$ and $\delta ^{\mu }=\delta u^{\mu
}=(d^{i},\delta ^{a})$ which emphasize that operators (\ref{dder}) and (%
\ref{ddif}) define, correspondingly, certain ``N--elongated'' partial derivatives and
differentials which are more
convenient for calculations on such nonholonomic manifolds.

Any N--connection also defines a N--connection curvature
\begin{equation*}
\mathbf{\Omega }=\frac{1}{2}\Omega _{ij}^{a}d^{i}\wedge d^{j}\otimes
\partial _{a},
\end{equation*}%
with N--connection curvature coefficients%
\begin{equation}
\Omega _{ij}^{a}=\delta _{\lbrack j}N_{i]}^{a}=\delta _{j}N_{i}^{a}-\delta _{i}N_{j}^{a}=
\frac{\partial N_{i}^{a}}{%
\partial x^{j}}-\frac{\partial N_{j}^{a}}{\partial x^{i}}+N_{i}^{b}\frac{%
\partial N_{j}^{a}}{\partial y^{b}}-N_{j}^{b}\frac{\partial N_{i}^{a}}{%
\partial y^{b}}.  \label{ncurv}
\end{equation}

The vielbeins (\ref{ddif}) satisfy the nonholonomy (equivalently, anholonomy) relations
\begin{equation}
\lbrack \mathbf{e}_{\alpha },\mathbf{e}_{\beta }]=\mathbf{e}_{\alpha }%
\mathbf{e}_{\beta }-\mathbf{e}_{\beta }\mathbf{e}_{\alpha }=W_{\alpha \beta }^{\gamma }\mathbf{e}_{\gamma }
\label{anhrel}
\end{equation}%
with (antisymmetric) nontrivial anholonomy coefficients $W_{ia}^{b}=\partial _{a}N_{i}^{b}$ and
$W_{ji}^{a}=\Omega _{ij}^{a}.$

\begin{definition}
A manifold \ $\mathbf{V}$ is called N--anholonomic if on the tangent space $T%
\mathbf{V}$ it is defined a local (nonintegrable) distribution (\ref{whitney}%
), i.e., $T\mathbf{V}$ is enabled with a N--connection and related nonholonomic vielbein structure
(\ref{anhrel}).
\end{definition}

We note that in this work we use boldfaced symbols for the spaces and geometric objects provided/adapted to a
N--connection structure. For
instance, a vector field $\mathbf{X}\in T\mathbf{V}$ \ is expressed $\mathbf{%
X}=(X,\ ^{\vee }X),$ or $\mathbf{X}=X^{\alpha }\mathbf{e}_{\alpha
}=X^{i}e_{i}+X^{a}e_{a},$ where $X=X^{i}e_{i}$ and $^{\vee
}X=X^{a}e_{a}$ state, respectively, the irreducible (adapted to
the N--connection structure) horizontal (h) and vertical (v)
components of the vector (which following refereces$^{6,7}$ is
called a distinguished vectors, in brief, d--vector). In a similar
fashion, the geometric objects on $\mathbf{V}$ like tensors,
spinors, connections, ... are called respectively d--tensors,
d--spinors, d--connections if they are adapted to the
N--connection splitting.

In the next two subsections we show how certain type of N--connection geometries can be
naturally derived from
Lagrange--Finsler geometry and in gravity theories.

\subsection{N--connections and Lagrangians}

We outline the main results on N--connections and almost
symplectic structures induced by regular Lagrangians$^{6-8}$. In
this case the N--anholono\-mic manifold $\mathbf{V}$ is to be
modelled on the tangent bundle $(TM, \pi, M),$ where $M$ is a
$n$--dimensional base manifold, $\pi$ is a surjective projection
and $TM$ is the total space. One denotes by
$\widetilde{TM}=TM\backslash \{0\}$ where $\{0\}$ means the null
section of map $\pi .$

A differentiable Lagrangian $L(x,y),$ i. e. a fundamental Lagrange function,
is defined by a map $L:(x,y)\in TM\rightarrow L(x,y)\in \mathbb{R}$ of class $%
\mathcal{C}^{\infty }$ on $\widetilde{TM}$ and continuous on the null section $0:M\rightarrow TM$ of $\pi .$ A
regular Lagrangian is with nondegenerated Hessian,
\begin{equation}
\ ^{(L)}g_{ij}(x,y)=\frac{1}{2}\frac{\partial ^{2}L(x,y)}{\partial y^{i}\partial y^{j}}  \label{lqf}
\end{equation}%
when $rank\left| g_{ij}\right| =n$ on $\widetilde{TM}.$

\begin{definition}
A Lagrange space is a pair $L^{n}=\left[ M,L(x,y)\right] $ with $\ ^{(L)}g_{ij}(x,y)$ being of constant
signature over $\widetilde{TM}.$
\end{definition}

The notion of Lagrange space was introduced by J. Kern$^{52}$ and
elaborated in details by the R. Miron's school on Finsler and
Lagrange geometry, see references$^{6,7}$, as a natural extension
of Finsler geometry$^{33,53-59}$ (see also references$^{45,46}$,
on Lagrange--Finsler supergeometry).

By straightforward calculations, there where proved the results:

\begin{enumerate}
\item The Euler--Lagrange equations%
\begin{equation*}
\frac{d}{d\tau }\left( \frac{\partial L}{\partial y^{i}}\right) -\frac{%
\partial L}{\partial x^{i}}=0
\end{equation*}%
where $y^{i}=\frac{dx^{i}}{d\tau }$ for $x^{i}(\tau )$ depending on parameter $\tau ,$ are equivalent to the
``nonlinear'' geodesic equations
\begin{equation*}
\frac{d^{2}x^{i}}{d\tau ^{2}}+2G^{i}(x^{k},\frac{dx^{j}}{d\tau })=0
\end{equation*}%
defining paths of a canonical semispray%
\begin{equation*}
S=y^{i}\frac{\partial }{\partial x^{i}}-2G^{i}(x,y)\frac{\partial }{\partial y^{i}}
\end{equation*}%
where
\begin{equation*}
2G^{i}(x,y)=\frac{1}{2}\ ^{(L)}g^{ij}\left( \frac{\partial ^{2}L}{\partial y^{i}\partial
x^{k}}y^{k}-\frac{\partial L}{\partial x^{i}}\right)
\end{equation*}%
with $^{(L)}g^{ij}$ being inverse to (\ref{lqf}).

\item There exists on $\widetilde{TM}$ a canonical N--connection $\ $%
\begin{equation}
^{(L)}N_{j}^{i}=\frac{\partial G^{i}(x,y)}{\partial y^{i}}  \label{cncl}
\end{equation}%
defined by the fundamental Lagrange function $L(x,y),$ which prescribes
nonholonomic frame structures of type (\ref{dder}) and (\ref{ddif}), $^{(L)}%
\mathbf{e}_{\nu }=(e_{i},^{\vee }e_{k})$ and $^{(L)}\mathbf{\vartheta }^{\mu
}=(\vartheta ^{i},^{\vee }\vartheta ^{k}).$ \footnote{%
On the tangent bundle the indices related to the base space run the same values as those related to fibers: we
can use the same symbols but have to distinguish like $^{\vee }e_{k}$ certain irreducible v--components with
respect to, (or for) N--adapted bases and cobases.}

\item The canonical N--connection (\ref{cncl}), defining $\ ^{\vee }e_{i},$
induces naturally an almost complex structure $\mathbf{F}:\chi (\widetilde{TM%
})\rightarrow \chi (\widetilde{TM}),$ where $\chi(\widetilde{TM})$ denotes the module of
vector fields on $\widetilde{TM},$%
\begin{equation*}
\mathbf{F}(e_{i})=\ ^{\vee }e_{i}\mbox{ and }\mathbf{F}(\ ^{\vee }e_{i})=-e_{i},
\end{equation*}%
when
\begin{equation}
\mathbf{F}=\ ^{\vee }e_{i}\otimes \vartheta ^{i}-e_{i}\otimes \ ^{\vee }\vartheta ^{i}
\label{acs1}
\end{equation}%
satisfies the condition $\mathbf{F\rfloor \ F=-I,}$ i. e. $F_{\ \ \beta
}^{\alpha }F_{\ \ \gamma }^{\beta }=-\delta _{\gamma }^{\alpha },$ where $%
\delta _{\gamma }^{\alpha }$ is the Kronecker symbol and ``$\mathbf{\rfloor } $"
denotes the interior product.

\item On $\widetilde{TM},$ there is a canonical metric structure%
\begin{equation}
\ ^{(L)}\mathbf{g}=\ ^{(L)}g_{ij}(x,y)\ \vartheta ^{i}\otimes \vartheta ^{j}+\ ^{(L)}g_{ij}(x,y)\ ^{\vee
}\vartheta ^{i}\otimes \ ^{\vee }\vartheta ^{j}  \label{slm}
\end{equation}%
constructed as a Sasaki type lift from $M.$
\end{enumerate}

One holds$^{6-8}$ the following

\begin{theorem}
$\label{slf1t}$The space $\left( \widetilde{TM},\mathbf{F,}^{(L)}\mathbf{g}%
\right) $ with almost complex form $\mathbf{F}$ (\ref{acs1}) defined by $%
^{(L)}N_{j}^{i}$, see (\ref{cncl}), and \ canonical metric structure $^{(L)}%
\mathbf{g}$ (\ref{slm}) is an almost K\"{a}hler space with almost symplectic
structure%
\begin{eqnarray}
^{(L)}\mathbf{\theta } &=&\ ^{(L)}\theta _{\alpha \beta }(x,y)\mathbf{%
\vartheta }^{\alpha }\wedge \mathbf{\vartheta }^{\beta }  \label{slf1} \\
&=&\ ^{(L)}g_{ij}(x,y)\ ^{\vee }\vartheta ^{i}\wedge \vartheta ^{j}.  \notag
\end{eqnarray}
\end{theorem}

\begin{proof}
It is evident if we define%
\begin{equation*}
^{(L)}\mathbf{\theta (X,Y)}\doteqdot \ \ ^{(L)}\mathbf{g(FX,Y)}
\end{equation*}%
and put $\mathbf{X=e}_{\alpha }$\textbf{\ }and $\mathbf{Y=e}_{\beta }.$ $\blacksquare$
\end{proof}
\vspace{3mm}

We conclude that any regular Lagrange mechanics can be geometrized as an almost K\"{a}hler space with
N--connection distribution. In a such Lagrange--K\"{a}hler nonholonomic manifold, the fundamental geometric
structures
(semispray, N--connection, almost complex structure and canonical metric on $%
\widetilde{TM})$ are defined by the fundamental Lagrange function $L(x,y).$

\begin{remark}
For applications in optics of nonhomogeneous media and gravity (see, for
instance, references$^{7,9,11,12}$) one considers metric forms of type $%
g_{ij}\sim e^{\lambda (x,y)}\ ^{(L)}g_{ij}(x,y) $ which can not be
derived from a mechanical Lagrangian. In the so--called
generalized Lagrange geometry one considers Sasaki type metrics
(\ref{slm}) with certain general coefficients both for the metric
and N--connection, i.e., when $ ^{(L)}g_{ij}\rightarrow $
$g_{ij}(x,y),$ and $^{(L)}N_{j}^{i}\rightarrow N_{j}^{i}(x,y).$
\footnote{ In this case, we can similarly define an almost
K\"{a}hler N--anholonomic space  $\left(
\widetilde{TM},\mathbf{F,\theta }\right) $ with the geometric
structures induced naturally by the N--connection.}
\end{remark}

\begin{remark}
Finsler geometry with the fundamental Finsler function $F(x,y),$ being homogeneous of type $F(x,\lambda
y)=\lambda F(x,y),$ for nonzero $\lambda
\in \mathbb{R},$ may be considered as a particular case of Lagrange geometry when $%
L=F^{2}.$ \footnote{%
In another turn, there is a proof$^{59}$ that any Lagrange
fundamental function $L$ can be modelled as a singular case in a
certain class Finsler geometries of extra dimension.} We shall
apply the methods of Finsler geometry and its almost K\"{a}hler
models in this work. Nevertheless, because the generalized
Lagrange spaces are very general ones enabled with N--anholonomic
structure inducing a corresponding almost symplectic structure we
shall emphasize just such geometric configurations.
\end{remark}

\begin{remark}
It is also proved that both generalized Lagrange and Finsler
geometries can be modelled on Riemannian--Cartan N--anholonomic
mani\-folds $^{13,30-32}$ if off--diagonal metrics and
N--connections are introduced into consideration.
\end{remark}

Now we shall demonstrate how N--anholonomic configurations can defined in gravity theories. In this case, it is
convenient to work on a general manifold $\mathbf{V}, \dim \mathbf{V}=n+m$ with global splitting, instead of
the tangent bundle $\widetilde{TM}.$

\subsection{N--connections in gravity}

Let us consider a metric structure on $\mathbf{V}$ with the coefficients defined with respect to a local
coordinate basis $du^{\alpha }=\left(
dx^{i},dy^{a}\right) ,$%
\begin{equation*}
\mathbf{g}=\underline{g}_{\alpha \beta }(u)du^{\alpha }\otimes du^{\beta }
\end{equation*}%
with
\begin{equation}
\underline{g}_{\alpha \beta }=\left[
\begin{array}{cc}
g_{ij}+N_{i}^{a}N_{j}^{b}h_{ab} & N_{j}^{e}h_{ae} \\
N_{i}^{e}h_{be} & h_{ab}%
\end{array}%
\right] .  \label{ansatz}
\end{equation}

A metric, for instance, parametrized in the form (\ref{ansatz})\ is generic off--diagonal if it can not be
diagonalized by any coordinate transforms. Performing a frame transform with the coefficients

\begin{eqnarray}
\mathbf{e}_{\alpha }^{\ \underline{\alpha }}(u) &=&\left[
\begin{array}{cc}
e_{i}^{\ \underline{i}}(u) & N_{i}^{b}(u)e_{b}^{\ \underline{a}}(u) \\
0 & e_{a}^{\ \underline{a}}(u)%
\end{array}%
\right] ,  \label{vt1} \\
\mathbf{e}_{\ \underline{\beta }}^{\beta }(u) &=&\left[
\begin{array}{cc}
e_{\ \underline{i}}^{i\ }(u) & -N_{k}^{b}(u)e_{\ \underline{i}}^{k\ }(u) \\
0 & e_{\ \underline{a}}^{a\ }(u)%
\end{array}%
\right] ,  \label{vt2}
\end{eqnarray}%
we write equivalently the metric in the form
\begin{equation}
\mathbf{g}=\mathbf{g}_{\alpha \beta }\left( u\right) \mathbf{\vartheta }%
^{\alpha }\otimes \mathbf{\vartheta }^{\beta }=g_{ij}\left( u\right) \vartheta ^{i}\otimes \vartheta
^{j}+h_{ab}\left( u\right) \ ^{\vee }\vartheta ^{a}\otimes \ ^{\vee }\vartheta ^{b},  \label{metr}
\end{equation}%
where $g_{ij}\doteqdot \mathbf{g}\left( e_{i},e_{j}\right) $ and $%
h_{ab}\doteqdot \mathbf{g}\left( e_{a},e_{b}\right) $ \ and
\begin{equation*}
\mathbf{e}_{\alpha }=\mathbf{e}_{\alpha }^{\ \underline{\alpha }}\partial _{%
\underline{\alpha }}\mbox{ and }\mathbf{\vartheta }_{\ }^{\beta }=\mathbf{e}%
_{\ \underline{\beta }}^{\beta }du^{\underline{\beta }}.
\end{equation*}%
are vielbeins of type (\ref{dder}) and (\ref{ddif}) defined for arbitrary $%
N_{i}^{b}(u).$ We can consider a special class of manifolds provided with a global splitting into conventional
``horizontal" and ``vertical" subspaces (\ref{whitney}) induced by the ``off--diagonal" terms $N_{i}^{b}(u)$ and
prescribed type of nonholonomic frame structure.

If the manifold $\mathbf{V}$ is (pseudo) Riemannian, there is a unique
linear connection (the Levi--Civita connection) $\nabla ,$ which is metric, $%
\nabla \mathbf{g=0,}$ and torsionless, $\ ^{\nabla }T=0.$ Nevertheless, the connection $\nabla $ is not adapted
to the nonintegrable distribution
induced by $N_{i}^{b}(u).$ In this case, \footnote{%
For instance, in order to construct exact solutions parametrized by generic off--diagonal metrics, or for
investigating nonholonomic frame structures in gravity models with nontrivial torsion.} it is more convenient to
work with more general classes of linear connections which are N--adapted but contain nontrivial torsion
coefficients because of nontrivial nonholonomy coefficients $W_{\alpha \beta }^{\gamma }$ (\ref{anhrel}).

For a splitting of a (pseudo) Riemannian--Cartan space of
dimension $(n+m)$ (under certain constraints, we can consider
(pseudo) Riemannian configurations), the Lagrange and Finsler type
geometries were modelled by N--anholonomic structures as exact
solutions of gravitational field equations $^{9-13,31,32}$. In
this paper, we shall concentrate on N--anholonomic almost complex
structures of vacuum gravity which can be naturally defined as
$(n+n)$ configurations, in general, embedded in certain spaces of
dimension $(n+m),$ $m\geq n$.

\section{Connections on Almost Symplectic N--an\-ho\-lonomic Manifolds}

The geometric constructions can be adapted to the N--connection structure:

\begin{definition}
A distinguished connection (d--connection) $\mathbf{D}$ on a manifold $\mathbf{V}$ is a linear connection conserving under parallelism the Whitney sum (\ref{whitney}) defining a
general N--connection. Equivalently, $\mathbf{D}\mathbf{P}=0$, $\mathbf{P}$ being the almost product structure
defined by the N--connection.
\end{definition}

The N--adapted components $\mathbf{\Gamma }_{\beta \gamma }^{\alpha }$ of a d-connection $\mathbf{D}_{\alpha
}=(\delta _{\alpha }\rfloor \mathbf{D})$
are defined by the equations%
\begin{equation*}
\mathbf{D}_{\alpha }\delta _{\beta }=\mathbf{\Gamma }_{\ \alpha \beta }^{\gamma }\delta _{\gamma },
\end{equation*}%
or
\begin{equation}
\mathbf{\Gamma }_{\ \alpha \beta }^{\gamma }\left( u\right) =\left( \mathbf{D%
}_{\alpha }\delta _{\beta }\right) \rfloor \delta ^{\gamma }.  \label{dcon1}
\end{equation}%
In its turn, this defines a N--adapted splitting into h-- and v--covariant derivatives,
$\mathbf{D}=D+\ ^{\vee
}D,$ where $D_{k}=\left( L_{jk}^{i},L_{bk\;}^{a}\right) $ and $\ ^{\vee }D_{c}=\left(
C_{jk}^{i},C_{bc}^{a}\right) $ are introduced as corresponding h- and
v--parametrizations of (\ref{dcon1}),%
\begin{equation*}
L_{jk}^{i}=\left( \mathbf{D}_{k}e_{j}\right) \rfloor \vartheta ^{i},\quad L_{bk}^{a}=\left(
\mathbf{D}_{k}e_{b}\right) \rfloor \vartheta ^{a},~C_{jc}^{i}=\left( \mathbf{D}_{c}e_{j}\right) \rfloor
\vartheta ^{i},\quad C_{bc}^{a}=\left( \mathbf{D}_{c}e_{b}\right) \rfloor \vartheta ^{a}.
\end{equation*}%
The components $\mathbf{\Gamma }_{\ \alpha \beta }^{\gamma }=\left(
L_{jk}^{i},L_{bk}^{a},C_{jc}^{i},C_{bc}^{a}\right) $ completely define a d--connection $\mathbf{D}$ on a
N--anholonomic manifold $\mathbf{V}.$

The simplest way to perform computations with d--connections is to use
N--adapted differential forms like $\mathbf{\Gamma }_{\beta }^{\alpha }=%
\mathbf{\Gamma }_{\beta \gamma }^{\alpha }\mathbf{\vartheta }^{\gamma }$ with the coefficients defined with
respect to (\ref{ddif}) and (\ref{dder}).

We shall say that a d--connection $\mathbf{D}$ preserves an almost symplectic 2--form,
 of Lagrange type $^{(L)}\mathbf{\theta }$ (\ref{slf1})
 (or any general one, $\mathbf{\theta }$) defined from a generalized
Lagrange
geometry or N--anholonomic gravity model, if%
\begin{equation}
\mathbf{D\theta =0}  \label{cond1s}
\end{equation}%
or
\begin{equation*}
\mathbf{Z}(\mathbf{\theta (X,Y)})=\mathbf{\theta (D}_{Z}\mathbf{X,Y)+\theta (X,\mathbf{D}_{Z}Y)}
\end{equation*}%
for any d--vector fields $\mathbf{X,Y,Z\in }T\mathbf{V.}$

\begin{theorem}
The torsion $\mathcal{T}^{\alpha }\doteqdot \mathbf{D\vartheta }^{\alpha }=d%
\mathbf{\vartheta }^{\alpha }+\Gamma _{\beta }^{\alpha }\wedge \mathbf{%
\vartheta }^{\beta }$ of a d--connection has the irreducible h- v-- components (d--torsions) with N--adapted
coefficients
\begin{eqnarray}
T_{\ jk}^{i} &=&L_{\ [jk]}^{i},\ T_{\ ja}^{i}=-T_{\ aj}^{i}=C_{\ ja}^{i},\
T_{\ ji}^{a}=\Omega _{\ ji}^{a},\   \notag \\
T_{\ bi}^{a} &=&T_{\ ib}^{a}=\frac{\partial N_{i}^{a}}{\partial y^{b}}-L_{\ bi}^{a},\ T_{\ bc}^{a}=C_{\
[bc]}^{a},  \label{dtors}
\end{eqnarray}

\noindent where $L_{\ [jk]}^{i}=L_{\ jk}^{i}-L_{\ kj}^{i}$ and so on.
\end{theorem}

\begin{proof}
By a straightforward calculation we can verify the formulas.$\blacksquare$
\end{proof}

\begin{remark}
\label{Levi-Civita}
The Levi--Civita linear connection $\nabla =\{^{\nabla }\mathbf{\Gamma }%
_{\beta \gamma }^{\alpha }\},$ with vanishing both torsion and
nonmetricity, is not adapted to the global splitting
(\ref{whitney}). In fact, if $\nabla$ was adapted, then $\nabla
\mathbf{P}=0$, $\mathbf{P}$ being the almost product structure
defined by the N--connection, and then, as $\nabla$ is
torsionless,  one obtains by means of the Lemma 2.1.6 of $^{60}$
that the Nijenhuis tensor field $N_{\mathbf{P}}$ vanishes, thus
proving that both vertical and horizontal distributions are
involutive in the sense of Frobenius theorem, which is not our
case of anholonomic manifolds. Then, we must look for another
connection to study the geometry of these manifolds.
\end{remark}

One holds:

\begin{proposition}
There is a preferred, canonical d--connection structure,$\ \widehat{\mathbf{D%
}}\mathbf{,}$ $\ $on N--aholonomic manifold $\mathbf{V}$ constructed only
from the metric and N--con\-nec\-ti\-on coefficients $%
[g_{ij},h_{ab},N_{i}^{a}]$ and satisfying the conditions $\widehat{\mathbf{D}%
}\mathbf{g}=0$ and $\widehat{T}_{\ jk}^{i}=0$ and $\widehat{T}_{\ bc}^{a}=0.$
\end{proposition}

\begin{proof}
By straightforward calculations with respect to the N--adapted bases (\ref%
{ddif}) and (\ref{dder}), we can verify that the connection
\begin{equation}
\widehat{\mathbf{\Gamma }}_{\beta \gamma }^{\alpha }=\ ^{\nabla }\mathbf{%
\Gamma }_{\beta \gamma }^{\alpha }+\ \widehat{\mathbf{P}}_{\beta \gamma }^{\alpha }  \label{cdc}
\end{equation}%
with the deformation d--tensor \footnote{$ \widehat{\mathbf{P}}_{\beta \gamma }^{\alpha }$ is a tensor field of
type (1,2). As is well known, the sum of a linear connection and  a tensor field of type (1,2) is a new linear
connection. }
\begin{equation*}
\widehat{\mathbf{P}}_{\beta \gamma }^{\alpha }=(P_{jk}^{i}=0,P_{bk}^{a}=%
\frac{\partial N_{k}^{a}}{\partial y^{b}},P_{jc}^{i}=-\frac{1}{2}%
g^{ik}\Omega _{\ kj}^{a}h_{ca},P_{bc}^{a}=0)
\end{equation*}%
satisfies the conditions of this Proposition. It should be noted that, in general, the components
$\widehat{T}_{\ ja}^{i},\ \widehat{T}_{\ ji}^{a}$ and $\widehat{T}_{\ bi}^{a}$ are not zero. This is an
anholonomic frame (or, equivalently, off--diagonal metric) effect.$\blacksquare$
\end{proof}
\vspace{3mm}

With respect to the N--adapted frames, the coefficients\newline
$\widehat{\mathbf{\Gamma }}_{\ \alpha \beta }^{\gamma }=\left( \widehat{L}%
_{jk}^{i},\widehat{L}_{bk}^{a},\widehat{C}_{jc}^{i},\widehat{C}%
_{bc}^{a}\right) $ are computed:
\begin{eqnarray}
\widehat{L}_{jk}^{i} &=&\frac{1}{2}g^{ir}\left( \frac{\delta g_{jr}}{%
\partial x^{k}}+\frac{\delta g_{kr}}{\partial x^{j}}-\frac{\delta g_{jk}}{%
\partial x^{r}}\right) ,  \label{candcon} \\
\widehat{L}_{bk}^{a} &=&\frac{\partial N_{k}^{a}}{\partial y^{b}}+\frac{1}{2}%
h^{ac}\left( \frac{\delta h_{bc}}{\partial x^{k}}-\frac{\partial N_{k}^{d}}{%
\partial y^{b}}h_{dc}-\frac{\partial N_{k}^{d}}{\partial y^{c}}h_{db}\right)
,  \notag \\
\widehat{C}_{jc}^{i} &=&\frac{1}{2}g^{ik}\frac{\partial g_{jk}}{\partial
y^{c}},  \notag \\
\widehat{C}_{bc}^{a} &=&\frac{1}{2}h^{ad}\left( \frac{\partial h_{bd}}{%
\partial y^{c}}+\frac{\partial h_{cd}}{\partial y^{b}}-\frac{\partial h_{bc}%
}{\partial y^{d}}\right) .  \notag
\end{eqnarray}%
For the canonical d--connection there are satisfied the conditions of
vanishing of torsion on the h--subspace and v--subspace, i.e., $\widehat{T}%
_{jk}^{i}=\widehat{T}_{bc}^{a}=0.$ In more general cases, such components of
torsion are not zero, for instance, the metric d--connections of type $%
\mathbf{\Gamma }_{\ \alpha \beta }^{\gamma }=\left( \widehat{L}%
_{jk}^{i}+l_{jk}^{i}(u),\widehat{L}_{bk}^{a},\widehat{C}_{jc}^{i},\widehat{C}%
_{bc}^{a}+c_{bc}^{a}(u)\right) $ is also compatible with metric (\ref{metr}) and has nontrivial $T_{jk}^{i}$ and
$\widehat{T}_{bc}^{a}.$

Let us consider a special case with $\dim \mathbf{V=}n+n,\ h_{ab}\rightarrow g_{ij}$ and $N_{i}^{a}\rightarrow
N_{\ i}^{j}$ in (\ref{metr}) when a tangent bundle structure is locally modelled on $\mathbf{V.\,}$\ We denote a
such space by $\widetilde{\mathbf{V}}_{(n,n)}.$ One holds:

\begin{theorem}
The canonical d--connection $\widehat{\mathbf{D}}$ (\ref{candcon}) for a
local modelling of a $\widetilde{TM}$ space on $\widetilde{\mathbf{V}}%
_{(n,n)}$ is defined by $\widehat{\mathbf{\Gamma }}_{\ \alpha \beta
}^{\gamma }=(\widehat{L}_{jk}^{i},\widehat{C}_{jk}^{i})\ $with%
\begin{equation}
\widehat{L}_{jk}^{i}=\frac{1}{2}g^{ir}\left( \frac{\delta g_{jr}}{\partial x^{k}}+\frac{\delta g_{kr}}{\partial
x^{j}}-\frac{\delta g_{jk}}{\partial x^{r}}\right) ,\widehat{C}_{jk}^{i}=\frac{1}{2}g^{ir}\left( \frac{\partial
g_{jr}}{\partial x^{k}}+\frac{\partial g_{kr}}{\partial x^{j}}-\frac{%
\partial g_{jk}}{\partial x^{r}}\right) .  \label{candcon1}
\end{equation}%
This d--connection is almost Hermitian, i.e., it is compatible with the almost Hermitian structure
$(\mathbf{g,F}),$ when
\begin{equation}
\widehat{\mathbf{D}}\mathbf{\theta }=0\mbox{ and }\widehat{\mathbf{D}}%
\mathbf{F}=0  \label{cond1}
\end{equation}%
for a 2--form \footnote{In an intrinsic way, $\theta (X,Y)=\mathbf{g}(\mathbf{F}X,Y)$.}
\begin{equation}
\mathbf{\theta }=\ \theta _{\alpha \beta }(x,y)\mathbf{\vartheta }^{\alpha }\wedge \mathbf{\vartheta }^{\beta
}=g_{ij}(x,y)\ ^{\vee }\vartheta ^{i}\wedge \vartheta ^{j}.  \notag
\end{equation}
\end{theorem}

\begin{proof}
It is similar to that for the Theorem \ref{slf1t}.$\blacksquare$
\end{proof}
\vspace{3mm}

On almost symplectic manifolds, usually there are considered symmetric linear connections. In our case, we can
always define a symmetric d--connection by taking the symmetric part \footnote{In coordinate-free notation,
$\mathbf{S}_{X}Y=\frac{1}{2}(\mathbf{D}_{X}Y+\mathbf{D}_{Y}X+[X,Y])$.} of $\mathbf{\Gamma }_{\ \alpha \beta
}^{\gamma },$
\begin{equation}
\mathbf{S}_{\ \alpha \beta }^{\gamma }=\frac{1}{2}\left( \mathbf{\Gamma }_{\ \alpha \beta }^{\gamma
}+\mathbf{\Gamma }_{\ \beta \alpha }^{\gamma }\right) ,  \label{a1}
\end{equation}%
where $\mathbf{\Gamma }_{\ \alpha \beta }^{\gamma }=(\widehat{L}%
_{jk}^{i}+l_{jk}^{i}(u),\widehat{C}_{jk}^{i}+c_{jk}^{i}(u)).\ $ On a N--anholonomic manifold
$\widetilde{\mathbf{V}}_{(n,n)},$ an almost symplectic form $\mathbf{\theta }$ is not closed, i. e.
$d\mathbf{\theta \neq
0.}$ But it may be closed under the action of N--adapted derivatives (\ref%
{dder}) and differentials (\ref{ddif}) when
\begin{equation*}
\delta \mathbf{\theta =}\delta (\theta _{\alpha \beta }(x,y)\mathbf{%
\vartheta }^{\alpha }\wedge \mathbf{\vartheta }^{\beta })=0,
\end{equation*}%
which means that
\begin{equation}
\mathbf{e}_{\gamma }\theta _{\alpha \beta }+\mathbf{e}_{\alpha }\theta _{\gamma \beta }+\mathbf{e}_{\beta
}\theta _{\alpha \gamma }=0.  \label{a2}
\end{equation}%
The condition (\ref{cond1s}) written in N--adapted bases results in
\begin{equation*}
\mathbf{e}_{\gamma }\theta _{\alpha \beta }=\mathbf{\Gamma }_{\ \alpha \gamma \beta }-\mathbf{\Gamma }_{\ \beta
\gamma \alpha }
\end{equation*}%
for $\mathbf{\Gamma }_{\ \alpha \gamma \beta }\doteqdot \theta _{\alpha \tau }\mathbf{\Gamma }_{\ \gamma \beta
}^{\tau }.$

\begin{definition}
An almost symplectic 2--form $\mathbf{\theta }$ is N--symplectic if it satisfies the conditions (\ref{a2}).
\end{definition}

There is a relation between the set of all d--connections $\mathbf{D}$ for which $\mathbf{D\theta =0}$ for any
given $\mathbf{\theta }$ and $\mathbf{N}
$ and the set of all symmetric connections on $\widetilde{\mathbf{V}}%
_{(n,n)}.$ By straightforward calculations we can verify that
\begin{equation}
\mathbf{\Gamma }_{\ \alpha \gamma \beta }=\frac{1}{2}(\mathbf{e}_{\alpha
}\theta _{\gamma \beta }-\mathbf{e}_{\gamma }\theta _{\alpha \beta }-\mathbf{%
e}_{\beta }\theta _{\alpha \gamma })+(\mathbf{S}_{\alpha \gamma \beta }-%
\mathbf{S}_{\gamma \beta \alpha }+\mathbf{S}_{\beta \gamma \alpha }) \label{a3}
\end{equation}%
is inverse to (\ref{a1}), which for almost symplectic $\theta _{\alpha \beta } $ satisfying the conditions
(\ref{a2}) simplifies to
\begin{equation*}
\mathbf{\Gamma }_{\ \alpha \gamma \beta }=\mathbf{e}_{\alpha }\theta _{\gamma \beta }+(\mathbf{S}_{\alpha \gamma
\beta }-\mathbf{S}_{\gamma \beta \alpha }+\mathbf{S}_{\beta \gamma \alpha }).
\end{equation*}

On holonomic manifolds with trivial N--connection, the formulas (\ref{a2})
and (\ref{a3}) transform into those from reference$^{3}$ with $\mathbf{e}%
_{\alpha }\rightarrow \partial /\partial u^{\alpha }.$ We may
conclude that N--anholonomic transforms map symplectic forms in
almost symplectic ones but preserve the main symmetry properties
and compatibility with the linear connection structure if the
computations are performed with respect to N--adapted bases.

\section{Curvature of N--symplectic d--Con\-nec\-ti\-ons}

Let $\mathbf{V}$ (or $\widetilde{\mathbf{V}}_{(n,n)}$) be an N--anholonomic manifold provided with a metric
d--connection $\mathbf{\Gamma }_{\gamma }^{\alpha }.$

\begin{definition}
A Fedosov N--anholonomic manifold is defined by an almost symplectic d--connection and almost complex structure
induced by the N--con\-nec\-ti\-on.
\end{definition}

\begin{definition}
A Lagrange--Fedosov manifold is a Fedosov N--anholonomic manifold with the N--connection and almost complex
structure defined by the fundamental Lagrange function, see Theorem \ref{slf1t}.
\end{definition}

The curvature of a symplectic d--connection $\mathbf{D}$ is defined by the
usual formula%
\begin{equation*}
\mathbf{R}(\mathbf{X},\mathbf{Y})\mathbf{Z}\doteqdot \mathbf{D}_{X}\mathbf{D}%
_{Y}\mathbf{Z}-\mathbf{D}_{Y}\mathbf{D}_{X}\mathbf{Z-D}_{[X,X]}\mathbf{Z.}
\end{equation*}%
Because on N--anholonomic spaces the ``simplest" adapted to the
N--connec\-ti\-on  induced almost complex structures is defined by
the canonical d--connec\-ti\-on, it is convenient to use it as a
symplectic d--connection.

By straightforward calculations we prove:

\begin{theorem}
The curvature $\mathcal{R}_{\ \beta }^{\alpha }\doteqdot \mathbf{D\Gamma }%
_{\beta }^{\alpha }=d\mathbf{\Gamma }_{\beta }^{\alpha }-\mathbf{\Gamma }%
_{\beta }^{\gamma }\wedge \mathbf{\Gamma }_{\gamma }^{\alpha }$ of a d--connection $\mathbf{\Gamma }_{\gamma
}^{\alpha }$ has the irreducible h- v-- components (d--curvatures) of $\mathbf{R}_{\ \beta \gamma \delta
}^{\alpha }$,%
\begin{eqnarray}
R_{\ hjk}^{i} &=&e_{k}L_{\ hj}^{i}-e_{j}L_{\ hk}^{i}+L_{\ hj}^{m}L_{\
mk}^{i}-L_{\ hk}^{m}L_{\ mj}^{i}-C_{\ ha}^{i}\Omega _{\ kj}^{a},  \notag \\
R_{\ bjk}^{a} &=&e_{k}L_{\ bj}^{a}-e_{j}L_{\ bk}^{a}+L_{\ bj}^{c}L_{\
ck}^{a}-L_{\ bk}^{c}L_{\ cj}^{a}-C_{\ bc}^{a}\Omega _{\ kj}^{c},  \notag \\
R_{\ jka}^{i} &=&e_{a}L_{\ jk}^{i}-D_{k}C_{\ ja}^{i}+C_{\ jb}^{i}T_{\
ka}^{b},  \label{dcurv} \\
R_{\ bka}^{c} &=&e_{a}L_{\ bk}^{c}-D_{k}C_{\ ba}^{c}+C_{\ bd}^{c}T_{\
ka}^{c},  \notag \\
R_{\ jbc}^{i} &=&e_{c}C_{\ jb}^{i}-e_{b}C_{\ jc}^{i}+C_{\ jb}^{h}C_{\
hc}^{i}-C_{\ jc}^{h}C_{\ hb}^{i},  \notag \\
R_{\ bcd}^{a} &=&e_{d}C_{\ bc}^{a}-e_{c}C_{\ bd}^{a}+C_{\ bc}^{e}C_{\ ed}^{a}-C_{\ bd}^{e}C_{\ ec}^{a}.  \notag
\end{eqnarray}
\end{theorem}

\begin{remark}
For an N--anholonomic manifold $\widetilde{\mathbf{V}}_{(n,n)}$ provided with N--sym\-plet\-ic canonical
d--connection $\widehat{\mathbf{\Gamma }}_{\ \gamma \alpha \beta }=\theta _{\gamma \tau }\widehat{\mathbf{\Gamma
}}_{\
\alpha \beta }^{\tau },$ see (\ref{candcon1}), the d--curvatures (\ref{dcurv}%
) reduces to three irreducible components
\begin{eqnarray}
R_{\ hjk}^{i} &=&e_{k}L_{\ hj}^{i}-e_{j}L_{\ hk}^{i}+L_{\ hj}^{m}L_{\
mk}^{i}-L_{\ hk}^{m}L_{\ mj}^{i}-C_{\ ha}^{i}\Omega _{\ kj}^{a},  \notag \\
R_{\ jka}^{i} &=&e_{a}L_{\ jk}^{i}-D_{k}C_{\ ja}^{i}+C_{\ jb}^{i}T_{\
ka}^{b},  \label{dcurv1} \\
R_{\ bcd}^{a} &=&e_{d}C_{\ bc}^{a}-e_{c}C_{\ bd}^{a}+C_{\ bc}^{e}C_{\ ed}^{a}-C_{\ bd}^{e}C_{\ ec}^{a}  \notag
\end{eqnarray}%
where all indices $i,j,k \ldots$ and $a,b, \ldots$ run the same values
but label the components with respect to different
h-- or v--frames.
\end{remark}

The indices of the components of the curvature tensor are lowered as%
\begin{equation*}
\mathbf{R}_{\ \tau \beta \gamma \delta }=\mathbf{\theta }_{\tau \alpha }%
\mathbf{R}_{\ \beta \gamma \delta }^{\alpha }.
\end{equation*}%
For Lagrange--Fedosov manifolds, the 2--form $\mathbf{\theta }_{\tau \alpha }$
has the coefficients defined by the metric structure and Lagrangian, see (%
\ref{slf1}). In this case we can apply the canonical d--connection and the d--metric for definition of the
curvature of symplectic d--connections.

Contracting respectively the components of (\ref{dcurv}) and (\ref{dcurv1}) we prove:

\begin{corollary}
The Ricci d--tensor $\mathbf{R}_{\alpha \beta }\doteqdot \mathbf{R}_{\
\alpha \beta \tau }^{\tau }$ has the irreducible h- v--components%
\begin{equation}
R_{ij}\doteqdot R_{\ ijk}^{k},\ \ R_{ia}\doteqdot -R_{\ ika}^{k},\ R_{ai}\doteqdot R_{\ aib}^{b},\
R_{ab}\doteqdot R_{\ abc}^{c}, \label{dricci}
\end{equation}%
for a general N--holonomic manifold $\mathbf{V,}$ and
\begin{equation}
R_{ij}\doteqdot R_{\ ijk}^{k},\ \ R_{ia}\doteqdot -R_{\ ika}^{k},\ \ R_{ab}\doteqdot R_{\ abc}^{c},
\label{dricci1}
\end{equation}%
for an N--anholonomic manifold $\widetilde{\mathbf{V}}_{(n,n)}.$
\end{corollary}

\begin{corollary}
The scalar curvature of a d--connection is
\begin{eqnarray*}
\overleftarrow{\mathbf{R}} &\doteqdot &\mathbf{g}^{\alpha \beta }\mathbf{R}%
_{\alpha \beta }=g^{ij}R_{ij}+h^{ab}R_{ab},\mbox{ for }\mathbf{V;} \\
&=&2g^{ij}R_{ij},\mbox{ for }\widetilde{\mathbf{V}}_{(n,n)}.
\end{eqnarray*}
\end{corollary}

\begin{corollary}
The Einstein d--tensor is computed $\mathbf{G}_{\alpha \beta }=\mathbf{R}%
_{\alpha \beta }-\frac{1}{2}\mathbf{g}_{\alpha \beta }\overleftarrow{\mathbf{%
R}}.$
\end{corollary}

In modern gravity theories, one considers more general linear connections generated by deformations of type
$\mathbf{\Gamma }_{\beta \gamma }^{\alpha }=\widehat{\mathbf{\Gamma }}_{\beta \gamma }^{\alpha
}+\mathbf{P}_{\beta \gamma }^{\alpha }$. We can split all geometric objects into canonical and post-canonical
pieces which results in N--adapted geometric constructions. For instance,
\begin{equation}
\mathcal{R}_{\ \beta }^{\alpha }=\widehat{\mathcal{R}}_{\ \beta }^{\alpha }+%
\mathbf{D}\mathcal{P}_{\ \beta }^{\alpha }+\mathcal{P}_{\ \gamma }^{\alpha }\wedge \mathcal{P}_{\ \beta
}^{\gamma }  \label{deformcurv}
\end{equation}%
for $\mathcal{P}_{\beta }^{\alpha }=\mathbf{P}_{\beta \gamma }^{\alpha }%
\mathbf{\vartheta }^{\gamma }.$ This way, for almost complex geometries, the d--tensors (\ref{dcurv1}) and
(\ref{dricci1}) can be redefined just for symmetrized d--connections compatible with the almost complex
structure.

\section{Einstein Flat N--Anholonomic Manifolds}

In terms of differential forms, the vacuum Einstein equations are written
\begin{equation}
\eta _{\alpha \beta \gamma }\wedge \widehat{\mathcal{R}}_{\ }^{\beta \gamma }=0,  \label{eecdc2}
\end{equation}%
where, for the volume 4--form $\eta \doteqdot \ast 1$ with the Hodge operator ``$\ast $", $\eta _{\alpha
}\doteqdot \mathbf{e}_{\alpha }\rfloor \eta ,$ $\eta _{\alpha \beta }\doteqdot \mathbf{e}_{\beta }\rfloor \eta
_{\alpha },$ $\eta _{\alpha \beta \gamma }\doteqdot \mathbf{e}_{\gamma }\rfloor \eta _{\alpha \beta },...$ and
$\widehat{\mathcal{R}}_{\ }^{\beta
\gamma }$ is the curvature 2--form. The deformation of connection (\ref{cdc}%
) defines a deformation of the curvature tensor of type (\ref{deformcurv}) but with respect to the curvature of
the Levi--Civita connection, $\ ^{\nabla }\mathcal{R}_{\ }^{\beta \gamma }.$ The gravitational field equations
(\ref{eecdc2}) transforms into
\begin{equation}
\eta _{\alpha \beta \gamma }\wedge \ ^{\nabla }\mathcal{R}_{\ }^{\beta \gamma }+\eta _{\alpha \beta \gamma
}\wedge \ ^{\nabla }\mathcal{Z}_{\ }^{\beta \gamma }=0,  \label{eecdc3}
\end{equation}%
where $^{\nabla }\mathcal{Z}_{\ \ \gamma }^{\beta }=\nabla \mathcal{P}_{\ \ \gamma }^{\beta }+\mathcal{P}_{\ \
\alpha }^{\beta }\wedge \mathcal{P}_{\ \ \gamma }^{\alpha }.$

A subclass of solutions of the gravitational field equations for the canonical d--connection defines also
solutions of the Einstein equations for the Levi--Civita connection if and only if
\begin{equation}
\eta _{\alpha \beta \gamma }\wedge \ ^{\nabla }\mathcal{Z}_{\ }^{\beta \gamma }=0.  \label{einstc}
\end{equation}%
This property is very important for constructing exact solutions
in Einstein and string gravity, parametrized by generic
off--diagonal metrics and anholonomic frames with associated
N--connection structure (see reviews of results in
references$^{30,31}$ and$^{32}$).

\subsection{The ansatz for metric}

In this subsection we investigate a class of five dimensional vacuum Einstein solutions with nontrivial
associated N--connection and generic off--diago\-nal metric. We analyze the conditions when such solutions
reduce to four dimensions and posses almost complex structure.

Let us consider a five dimensional ansatz for the metric (\ref{metr}) and frame (\ref{ddif}) when $u^{\alpha
}=(x^{i},y^{4}=v,y^{5});i=1,2,3$ and the coefficients
\begin{eqnarray}
g_{ij} &=&diag[g_{1}=\pm \varpi (x^{k},v),\varpi
(x^{k},v)g_{2}(x^{2},x^{3}),\varpi (x^{k},v)g_{3}(x^{2},x^{3})],  \notag \\
h_{ab} &=&diag[\varpi (x^{k},v)h_{4}(x^{k},v),\varpi
(x^{k},v)h_{5}(x^{k},v)],  \notag \\
N_{i}^{4} &=&w_{i}(x^{k},v),N_{i}^{5}=n_{i}(x^{k},v)  \label{ansatz1}
\end{eqnarray}%
are some functions of necessary smooth class. The partial derivative are briefly denoted $a^{\bullet }=\partial
a/\partial x^{2},a^{^{\prime }}=\partial a/\partial x^{3},a^{\ast }=\partial a/\partial v.$

\begin{theorem}
\label{theq}The vacuum Einstein equations (\ref{eecdc2}) for the
ca\-no\-nical d--con\-nec\-ti\-on (\ref{cdc}) constructed from data (\ref%
{ansatz1}) are equivalent to the system of equations
\begin{eqnarray}
g_{3}^{\bullet \bullet }-\frac{g_{2}^{\bullet }g_{3}^{\bullet }}{2g_{2}}-%
\frac{(g_{3}^{\bullet })^{2}}{2g_{3}}+g_{2}^{^{\prime \prime }}-\frac{%
g_{2}^{\prime }g_{3}^{\prime }}{2g_{3}}-\frac{(g_{2}^{^{\prime }})^{2}}{%
2g_{2}} &=&0,  \label{ep1a} \\
h_{5}^{\ast \ast }-h_{5}^{\ast }(\ln \left| \sqrt{\left| h_{4}h_{5}\right| }%
\right| )^{\ast } &=&0,  \label{ep2a} \\
w_{i}\beta +\alpha _{i} &=&0,  \label{ep3a} \\
n_{i}^{\ast \ast }+n_{i}^{\ast } &=&0,  \label{ep4a}
\end{eqnarray}%
where%
\begin{eqnarray}
\alpha _{i} &=&\partial _{i}h_{5}^{\ast }-h_{5}^{\ast }\partial _{i}\ln \left| \sqrt{\left| h_{4}h_{5}\right|
}\right| ,\ \beta =h_{5}^{\ast \ast
}-h_{5}^{\ast }\left[ \ln \left| \sqrt{\left| h_{4}h_{5}\right| }\right| %
\right] ^{\ast },  \notag \\
\ \gamma &=&3h_{5}^{\ast }/2h_{5}-h_{4}^{\ast }/h_{4}  \label{abc}
\end{eqnarray}%
$h_{4}^{\ast }\neq 0$ and $h_{5}^{\ast }\neq 0$ and the functions $h_{4}$ and $\varpi $ must satisfy certain
additional conditions
\begin{equation}
\widehat{\delta }_{i}h_{4}=0\mbox{ and }\widehat{\delta }_{i}\varpi =0, \label{cond3}
\end{equation}%
for any $\zeta _{i}(x^{k},v)$ defining $\widehat{\delta }_{i}=\partial _{i}-(w_{i}+\zeta _{i})\partial
_{4}+n_{i}\partial _{5}.$
\end{theorem}

\begin{proof}
It is a straightforward calculation, see similar ones
in$^{31,9,11}$.$\blacksquare$
\end{proof}
\vspace{3mm}

We note that the conditions (\ref{cond3}) are satisfied if
\begin{equation}
\varpi ^{q_{1}/q_{2}}=h_{4}  \label{cond6}
\end{equation}%
for some nonzero integers $q_{1}$ and $q_{2}$ and $\zeta _{i}$ defined from
the equations%
\begin{equation}
\partial _{i}\varpi -(w_{i}+\zeta _{i})\varpi ^{\ast }=0.  \label{confeq}
\end{equation}

\begin{remark}
\label{rtheq} Under the conditions of the Theorem \ref{theq}, we can also
consider d--metrics with $h_{5}^{\ast }=0$ for such functions $%
h_{4}=h^{\#}(x^{i},v)$ when
\begin{equation*}
\lim_{h_{5}^{\ast }\rightarrow 0}\left\{ h_{5}^{\ast }[\ln \left| \sqrt{%
\left| h^{\#}h_{5}\right| }\right| ]^{\ast }\right\} \rightarrow 0
\end{equation*}%
and
\begin{equation*}
\lim_{h_{5}^{\ast }\rightarrow 0}\left\{ h_{5}^{\ast }\partial _{i}\ln \left| \sqrt{\left| h^{\#}h_{5}\right|
}\right| \right\} \rightarrow 0.
\end{equation*}%
In this cases, the equations (\ref{ep2a}) and (\ref{ep3a}) will be satisfied by any $h^{\#}(x^{i},v)$ and
$w_{i}(x^{i},v)$ and we may take $n_{i}^{\ast }=n_{[1]i}(x^{i})h^{\#}(x^{i},v)$ in order to satisfy
(\ref{ep4a}).
\end{remark}

\begin{theorem}
\label{texs}The system of gravitational field equations (\ref{eecdc2}) for the ansatz (\ref{ansatz1}) can be
solved in general form if there are given
certain values of functions $g_{2}(x^{2},x^{3})$ (or, inversely, $%
g_{3}(x^{2},x^{3})$), $h_{4}(x^{i},v)$ (or, inversely, $h_{5}(x^{i},v)$).
\end{theorem}

\begin{proof}
We outline the main steps of constructing exact solutions proving
this Theorem, see detailed computations presented in the Proof of
Theorem 4.3 from reference$^{31}$.

\begin{itemize}
\item The general solution of equation (\ref{ep1a}) can be written in the form
\begin{equation}
\lambda =g_{[0]}\exp [a_{2}\widetilde{x}^{2}\left( x^{2},x^{3}\right) +a_{3}%
\widetilde{x}^{3}\left( x^{2},x^{3}\right) ],  \label{solricci1a}
\end{equation}%
were $g_{[0]},a_{2}$ and $a_{3}$ are some constants and the functions $%
\widetilde{x}^{2,3}\left( x^{2},x^{3}\right) $ define any coordinate transforms $x^{2,3}\rightarrow
\widetilde{x}^{2,3}$ for which the two dimensional line element becomes conformally flat, i.e.,\begin{equation}
g_{2}(x^{2},x^{3})(dx^{2})^{2}+g_{3}(x^{2},x^{3})(dx^{3})^{2}\rightarrow
\lambda (x^{2},x^{3})\left[ (d\widetilde{x}^{2})^{2}+\epsilon (d\widetilde{x}%
^{3})^{2}\right] ,  \label{con10}
\end{equation}%
where $\epsilon =\pm 1$ for a corresponding signature. In coordinates $%
\widetilde{x}^{2,3},$ the equation (\ref{ep1a}) transform into%
\begin{equation*}
\lambda \left( \lambda ^{\bullet \bullet }+\lambda ^{\prime \prime }\right) -\lambda ^{\bullet }-\lambda
^{\prime }=0
\end{equation*}%
or%
\begin{equation}
\ddot{\psi}+\psi ^{\prime \prime }=0,  \label{auxeq01}
\end{equation}%
for $\psi =\ln |\lambda |.$ There are three alternative possibilities to
generate solutions of (\ref{ep1a}). For instance, we can prescribe that $%
g_{2}=g_{3}$ and get the equation (\ref{auxeq01}) for $\psi =\ln |g_{2}|=\ln
|g_{3}|.$ If we suppose that $g_{2}^{^{\prime }}=0,$ for a given $%
g_{2}(x^{2}),$ we obtain from (\ref{ep1a})%
\begin{equation*}
g_{3}^{\bullet \bullet }-\frac{g_{2}^{\bullet }g_{3}^{\bullet }}{2g_{2}}-%
\frac{(g_{3}^{\bullet })^{2}}{2g_{3}}=0
\end{equation*}%
which can be integrated exactly. Similarly, we can generate solutions for a prescribed $g_{3}(x^{3})$ in the
equation
\begin{equation}
g_{2}^{^{\prime \prime }}-\frac{g_{2}^{^{\prime }}g_{3}^{^{\prime }}}{2g_{3}}%
-\frac{(g_{2}^{^{\prime }})^{2}}{2g_{2}}=0.  \label{aux4}
\end{equation}%
We note that a transform (\ref{con10}) is always possible for 2D metrics and the explicit form of solutions
depends on chosen system of 2D coordinates
and on the signature $\epsilon =\pm 1.$ In the simplest case, the equation (%
\ref{ep1a}) is solved by arbitrary two functions $g_{2}(x^{3})$ and $%
g_{3}(x^{2}).$

\item The equation (\ref{ep2a}) relates two functions $h_{4}\left( x^{i},v\right) $ and $h_{5}\left(
x^{i},v\right) $ following two possibilities:

a) to compute
\begin{eqnarray}
\sqrt{|h_{5}|} &=&h_{5[1]}\left( x^{i}\right) +h_{5[2]}\left( x^{i}\right) \int \sqrt{|h_{4}\left(
x^{i},v\right) |}dv,~h_{4}^{\ast }\left(
x^{i},v\right) \neq 0;  \notag \\
&=&h_{5[1]}\left( x^{i}\right) +h_{5[2]}\left( x^{i}\right) v,\ h_{4}^{\ast }\left( x^{i},v\right) =0,
\label{p2}
\end{eqnarray}%
for some functions $h_{5[1,2]}\left( x^{i}\right) $ stated by boundary conditions;

b) or, inversely, to compute $h_{4}$ for a given $h_{5}\left( x^{i},v\right)
,h_{5}^{\ast }\neq 0,$%
\begin{equation}
\sqrt{|h_{4}|}=h_{[0]}\left( x^{i}\right) (\sqrt{|h_{5}\left( x^{i},v\right) |})^{\ast },  \label{p1}
\end{equation}%
with $h_{[0]}\left( x^{i}\right) $ given by boundary conditions.

\item The exact solutions of (\ref{ep3a}) for $\beta \neq 0$ are defined from an algebraic equation, $w_{i}\beta
+\alpha _{i}=0,$ where the
coefficients $\beta $ and $\alpha _{i}$ are computed as in formulas (\ref%
{abc}) by using the solutions for (\ref{ep1a}) and (\ref{ep2a}). The general solution is
\begin{equation}
w_{k}=\partial _{k}\ln [\sqrt{|h_{4}h_{5}|}/|h_{5}^{\ast }|]/\partial _{v}\ln [\sqrt{|h_{4}h_{5}|}/|h_{5}^{\ast
}|],  \label{w}
\end{equation}%
with $\partial _{v}=\partial /\partial v$ and $h_{5}^{\ast }\neq 0.$ If $%
h_{5}^{\ast }=0,$ or even $h_{5}^{\ast }\neq 0$ but $\beta =0,$ the coefficients $w_{k}$ could be arbitrary
functions on $\left( x^{i},v\right) . $ \ For the vacuum Einstein equations this is a degenerated case imposing
the compatibility conditions $\beta =\alpha _{i}=0,$ which are satisfied, for instance, if the $h_{4}$ and
$h_{5}$ are related as in the formula (\ref{p1}) but with $h_{[0]}\left( x^{i}\right) =const.$

\item Having defined $h_{4}$ and $h_{5}$ and computed $\gamma $ from (\ref%
{abc}) we can solve the equation (\ref{ep4a}) by integrating on variable ``$%
v $" the equation $n_{i}^{\ast \ast }+\gamma n_{i}^{\ast }=0.$ The exact solution is
\begin{eqnarray}
n_{k} &=&n_{k[1]}\left( x^{i}\right) +n_{k[2]}\left( x^{i}\right) \int
[h_{4}/(\sqrt{|h_{5}|})^{3}]dv,~h_{5}^{\ast }\neq 0;  \notag \\
&=&n_{k[1]}\left( x^{i}\right) +n_{k[2]}\left( x^{i}\right) \int
h_{4}dv,\qquad ~h_{5}^{\ast }=0;  \label{n} \\
&=&n_{k[1]}\left( x^{i}\right) +n_{k[2]}\left( x^{i}\right) \int [1/(\sqrt{%
|h_{5}|})^{3}]dv,~h_{4}^{\ast }=0,  \notag
\end{eqnarray}%
for some functions $n_{k[1,2]}\left( x^{i}\right) $ stated by boundary conditions.
\end{itemize}

The exact solution of (\ref{confeq}) is given by some functions $\zeta _{i}=\zeta _{i}\left( x^{i},v\right) $ if
\ both $\partial _{i}\varpi =0$ and $\varpi ^{\ast }=0,$ we chose $\zeta _{i}=0$ for $\varpi =const,$ and
\begin{eqnarray}
\zeta _{i} &=&-w_{i}+(\varpi ^{\ast })^{-1}\partial _{i}\varpi ,\quad \varpi
^{\ast }\neq 0,  \label{confsol} \\
&=&(\varpi ^{\ast })^{-1}\partial _{i}\varpi ,\quad \varpi ^{\ast }\neq 0,%
\mbox{ for vacuum solutions}.  \notag
\end{eqnarray}%
$\blacksquare$
\end{proof}

The Theorem \ref{texs} states a general method of constructing five dimensional exact solutions in various
gravity models with generic off--diagonal metrics, nonholonomic frames and, in general, with nontrivial torsion.
Such solutions are with associated N--connection structure. This method can be also applied in order to
generate, for instance, certain Finsler or Lagrange configurations as v-irreducible components, or for a certain
class of conformal factors $\varpi (x^{i},v)$ for both h-- and v--irreducible components. The five dimensional
ansatz can not be used to generate directly standard Finsler or Lagrange geometries because the dimension of
such spaces can not be an odd number. Nevertheless, the anholonomic frame method can be applied in order to
generate four dimensional exact solutions containing Finsler--Lagrange configurations. For
instance, a four dimensional configuration can be defined just by an ansatz (%
\ref{metr}) with the data (\ref{ansatz1}) where the coefficients do not depend on coordinate $x^{1}$ and the
metric is stated to be four dimensional with the conformal factor $\varpi (x^{2},x^{3},v).$

\subsection{An example of induced almost K\"{a}hler gravity}

Let us consider a four dimensional ansatz which may mimic under certain constraints a generalized Lagrange
geometry and induced almost K\"{a}hler
structure in Riemann--Cartan space:%
\begin{eqnarray*}
\mathbf{g} &=&\varpi (x^{2},x^{3},v)[g_{22}\left( x^{2},x^{3}\right)
dx^{2}\otimes dx^{2}+g_{33}\left( x^{2},x^{3}\right) dx^{3}\otimes dx^{3} \\
&&+h_{44}\left( x^{2},x^{3},v\right) \ \delta y^{4}\otimes \ \delta y^{4}+h_{55}\left( x^{2},x^{3},v\right) \
\delta y^{5}\otimes \ \delta y^{5}],
\end{eqnarray*}%
where%
\begin{eqnarray*}
\delta y^{4} &=&dv+w_{2}\left( x^{2},x^{3},v\right) dx^{2}+w_{3}\left(
x^{2},x^{3},v\right) dx^{3}, \\
\delta y^{5} &=&dy^{5}+n_{2}\left( x^{2},x^{3},v\right) dx^{2}+n_{3}\left( x^{2},x^{3},v\right) dx^{3}.
\end{eqnarray*}%
This d--metric will define a class of vacuum solutions of the Einstein
equations if the coefficients are subjected to the conditions of the Theorem%
\ref{texs}, when the dependence on coordinate $x^{1}$ is eliminated. We put $%
g_{22}=g(x^{3})$ and $g_{33}=0$ to be a solution of (\ref{ep1a}) in the form (\ref{aux4}),
i.e.,\begin{equation*}
2gg^{^{\prime \prime }}-(g^{^{\prime }})^{2}=0
\end{equation*}%
and choose $h_{5}=0$ and
\begin{equation*}
h_{4}=h^{\#}(x^{3},v)=\frac{a^{2}}{\left| g(x^{3})\times v\right| }g(x^{3})
\end{equation*}%
for $a=const,$ which satisfies (\ref{ep2a}), see Remark \ref{rtheq}. Taking any functions
$w_{2,3}(x^{2},x^{3},v)$ and $n_{2,3}(x^{2},x^{3},v)$ satisfying
\begin{equation*}
n_{2,3}^{\ast }=n_{2,3[0]}(x^{2},x^{3})h^{\#}(x^{3},v)
\end{equation*}%
we solve respectively the equations (\ref{ep3a}) and (\ref{ep4a}). We may take
\begin{equation*}
\varpi =\varpi ^{\#}(x^{3},v)=\left[ h^{\#}(x^{3},v)\right] ^{q_{2}/q_{1}}
\end{equation*}%
like for (\ref{confsol}). All such functions define a vacuum Einstein
d--metric%
\begin{equation*}
\mathbf{g}=\varpi ^{\#}(x^{3},v)[g\left( x^{3}\right) dx^{2}\otimes dx^{2}+%
\frac{a^{2}}{\left| g(x^{3})\times v^{2}\right| }g(x^{3})\delta y^{4}\otimes \ \delta y^{4}],
\end{equation*}%
modelling an embedded generalized Lagrange geometry (it is a particular case of d--metrics
 considered in$^{8}$, see formula (6.3), which in our case is derived from a gravity model). We construct a conformal
almost
K\"{a}hler geometry if we consider%
\begin{equation*}
\mathbf{\theta }=\varpi ^{\#}(x^{3},v)g\left( x^{3}\right) \frac{a}{\sqrt{%
\left| g(x^{3})\times v^{2}\right| }}\delta y^{4}\wedge dx^{2}
\end{equation*}%
and
\begin{eqnarray*}
\mathbf{F} &\mathbf{=}&\frac{\sqrt{\left| g(x^{3})\times v^{2}\right| }}{a}%
\left( \frac{\partial }{\partial v}\otimes dx^{2}+\frac{\partial }{\partial
y^{5}}\otimes dx^{3}\right)  \\
\noalign{\smallskip}
&&-\frac{a}{\sqrt{\left| g(x^{3})\times v^{2}\right| }}\left( \frac{\delta }{%
\partial x^{2}}\otimes dv+\frac{\delta }{\partial x^{3}}\otimes
dy^{5}\right) .
\end{eqnarray*}%
Finally, we note that if we choose the functions
$w_{2,3}(x^{2},x^{3},v)$ and \newline $n_{2,3}(x^{2},x^{3},v)$ to
parametrize a noncommutative structure, this vacuum gravitational
space will possess a noncommutative symmetry like in$^{31,32}$. An
alternative class of solutions can be generated if we put certain
boundary conditions (for instance, for $v=t$ treated as a timelike
coordinate, and one of the space coordinates $x^{2},x^{3},y^{5}$
running to infinite) when the N--connection coefficients possess a
Lie algebra symmetry. In this case, we generate an explicit
example of vacuum gravitational fields (in general, with
nontrivial torsion) possessing Lie symmetries$^{61}$. We can
select such values of $w_{2,3}$ and $n_{2,3} $ when the conditions
(\ref{einstc}) are satisfied and the solutions coincide with those
for the Levi Civita connection, but this is a very restricted case
of N--connection geometry and associated almost complex
structures.

\subsection*{Acknowledgment}

F. E.'s research was partially supported by the Spanish Ministerio de Ciencia y Tecnolog\'{\i}a (BFM
2002-00141). R. S. is partially supported by ULE2003-02 (Spain). S. V. is grateful to the Vicerrectorado de
Investigaci\'{o}n de la Universidad de Cantabria for financial support
 and C. Tanasescu for support and kind hospitality.
\vspace{5mm}

\footnotesize
\noindent $^{1}$B. V. Fedosov, J. Diff. Geom. \textbf{40}, 213 (1994).\\
$^{2}$B. V. Fedosov, {\em Deformation Quantization and Index
Theory} (Akad. Verlag. Berlin, 1996).\\
$^{3}$I. Gelfand, V. Retakh and M. Shubin, Adv. Math.
\textbf{136},
104 (1998).\\
$^{4}$M. De Wilde and P. B. A. Lecompte, Ann. Inst. Fourier
\textbf{35}, 117 (1985).\\
$^{5}$C. Castro, J. Geom. Phys. {\bf 33}, 173 (2000).\\
$^{6}$R. Miron and M. Anastasiei, {\em Vector Bundles and Lagrange
Spaces with Applications to Relativity} (Geometry Balkan Press,
Bukharest, 1997, translation from Romanian of Etitura Academiei
Romane, 1984).\\
$^{7}$R. Miron and M. Anastasiei, {\em The Geometry of Lagrange
Spaces: Theory and Applications} (Kluwer, 1994).\\
$^{8}$R. Miron, D. Hrimiuc, H. Shimada and V. S. Sabau, {\em The
Geometry of Hamilton and Lagrange Spaces} (Kluwer, 2000).\\
$^{9}$S. Vacaru, JHEP {\bf 04}, 009 (2001).\\
$^{10}$S. Vacaru,  Ann. Phys. (N. Y.) \textbf{290}, 83 (2001).\\
$^{11}$S. Vacaru and O. Tintareanu-Mircea, Nucl. Phys. B
\textbf{626}, 239 (2002).\\
$^{12}$S. Vacaru and D. Singleton,  Class. Quant. Gravity,
\textbf{19}, 2793 (2002).\\
$^{13}$S. Vacaru and H. Dehnen,  Gen. Rel. Grav. \textbf{35}, 209
(2003).\\
$^{14}$I. Vaisman, {\em Lectures on the Geometry of Poisson
Manifolds} (Birkhauser, Verlag, Basel, 1994).\\
$^{15}$H. C. Lee,  Amer. J. Math. \textbf{65}, 433 (1943).\\
$^{16}$P. Libermann, Ann. Mat. Pura Appl. \textbf{36}, 27 (1954).\\
$^{17}$Yu. I. Levin, Dokl. Akad. Nauk SSSR. \textbf{128}, 668
(1959). [in Russian].\\
$^{18}$P. Tondeur, Comment. Math. Helv. \textbf{36}, 234 (1961).\\
$^{19}$I. Vaisman, Monatsh. Math. \textbf{100}, 299 (1985).\\
$^{20}$D. Leites, E. Poletaeva and V. Serganova, J. Nonlin. Math.
Phys. \textbf{9}, 394 (2002).\\
$^{21}$D. Leites, The Riemann tensor for nonholonomic manifolds,
(unpublished, math. RT/ 0202213).\\
$^{22}$A. M. Vershik, Classical and Nonclassical Dynamics with
Constraints, in: Yu. G. Borisovish and Yu. F. Gliklih (eds.),
Geometry and Topology in Global Nonliner Problems, Novoe Glob.
Anal. Voronezh Gos. Univ. Voronezh, 1984; Engl. trans. in:  Global
Analysis -- Studies and Applications. I. (A collection of articles
translated from the Russian), Lect. Notes in Mathematics {\bf
1108}, 278 (1984).\\
$^{23}$R. Miron, Part. II. An. Univ. Bucuresti, Math. Inform.
\textbf{50}, 93 (2001).\\
$^{24}$R. Miron, Algebras, Groups. Geom. \textbf{17}, 283 (2001).\\
$^{25}$R. Miron, An. \c Sti. Univ. ``Al. I. Cuza'' Ia\c si. Sect.
I (N. S.) \textbf{3}, 171 (1957).\\
$^{26}$R. Miron,  Acad. R. P. Rom\^{\i}ne. Fil. Ia\c si. Stud.
Cerc. \c Sti. Mat. \textbf{8}, 49 (1957).\\
$^{27}$R. Miron, An. \c Sti. Univ. ``Al. I. Cuza'' Ia\c si. Sec\c
t. I. (N. S.) \textbf{2}, 85 (1956).\\
$^{28}$V. Cruceanu, Acad. R. P. Rom\`{\i}ne Fil. Ia\c si Stud.
Cerc. \c Sti. Mat. \textbf{11}, 343 (1960).\\
$^{29}$F. Cantrijn, M. de Le\'{o}n, J. C. Marrero and D.
Mart\'{\i}n de Diego. Nonlinearity \textbf{13}, 1379 (2000).\\
$^{30}$S. Vacaru, Generalized Finsler Geometry in Einstein, String
and Metric--Affine Gravity, (unpublished, hep-th/0310132).\\
$^{31}$S. Vacaru, E. Gaburov and D. Gontsa, A Method of
Constructing Off--Diagonal Solutions in Metric--Affine and String
Gravity, (unpublished, hep--th/0310133).\\
$^{32}$S. Vacaru and E. Gaburov, Noncomnmutative Symmetries and
Stability of Black Ellipsoids in Metric--Affine and String
Gravity, (unpublished, hep--th/0310134).\\
$^{33}$E. Cartan, {\em Les Espaces de Finsler} (Herman, Paris,
1935).\\
$^{34}$A. Kawaguchi,  Akad. Wetensch. Amsterdam Proc. \textbf{40},
596 (1937).\\
$^{35}$A. Kawaguchi, On the Theory of Non--linear Connections I.
Tensor, N. S. {\bf 2}, 123 (1952).\\
$^{36}$ A. Kawaguchi, On the Theory of Non--linear Connections II.
Tensor, N. S. \textbf{6}, 165 (1956).\\
$^{37}$C. Ehresmann, Colloque de topologie (espaces fibr\'{e}s),
Bruxelles 1950, 29 (1951).\\
$^{38}$W. Barthel, J. Reine Angew. Math. \textbf{212}, 120 (1963).\\
$^{39}$M. de Le\'{o}n and C. Villaverde, C. R. Acad. Sci. Paris,
serie I, \textbf{293}, 51 (1981).\\
$^{40}$M. de Le\'{o}n, J. Mar\'{\i}n-Solano and J. Marrero, An.
Fis. Monogr. CIEMAT (Madrid) {\bf 2}, 73 (1995).\\
$^{41}$F. Etayo, Riv. Mat. Univ. Parma \textbf{17}, 131 (1991).\\
$^{42}$A. Ferr\'{a}ndez, Rev. Roumaine Math. Pures Appl.
\textbf{29} no. 3, 225 (1984).\\
$^{43}$S. Vacaru, J. Math. Phys. \textbf{37}, 508 (1996).\\
$^{44}$S. Vacaru and N. Vicol, Int. J. Math. and Math. Sciences
(IJMMS) \textbf{23}, 1189 (2004).\\
$^{45}$S. Vacaru, Nucl. Phys. B \textbf{434}, 590 (1997).\\
$^{46}$S. Vacaru, Interactions, {\em Strings and Isotopies in
Higher Order An\-iso\-tro\-pic Superspaces} (Hadronic Press, Palm
Harbor, FL, USA, 1998).\\
$^{47}$S. Vacaru, Phys. Lett. B \textbf{498}, 74 (2001).\\
$^{48}$J. Grifone, Structure presque-tangente et connexions, I
 and II, Ann. Inst. Fourier (Grenoble) \textbf {22}, no.1 287 and
 no. 3  291 (1972).\\
$^{49}$F. Etayo,  Rev. Acad. Canaria Cienc. \textbf{5}, 125
(1993).\\
$^{50}$K. Yano and S. Ishihara,  {\em Tangent and cotangent
bundles. Differential geometry}  (Marcel Dekker, Inc., New York,
1973).\\
$^{51}$M. Crampin, Math. Proc. Cambridge Philos. Soc. \textbf{94}
no. 1, 125 (1983).\\
$^{52}$J. Kern, Arch. Math. \textbf{25}, 438 (1974).\\
$^{53}$P. Finsler, \textit{\"{U}ber Kur\-ven und Fl\"{a}chen in
Allgemeiner R\"{a}men} (Dis\-ser\-ta\-ti\-on G\"{o}ttingen, 1918,
reprinted Basel: Birkh\"{a}user, 1951).\\
$^{54}$H. Rund, \textit{The Differential Geometry of Finsler
Spaces} (Berlin: Sprin\-ger--Verlag, 1959).\\
$^{55}$G. S. Asanov, \textit{Finsler Geometry, Relativity and
Gauge Theories} (Klu\-wer, Reidel, 1985).\\
$^{56}$M. Matsumoto, \ \textit{Foundations of Finsler Geometry and
Special Finsler Spaces} (Kaisisha: Shigaken, 1986).\\
$^{57}$A. Bejancu, \textit{Finsler Geometry and Applications}
(Ellis Horwood, Chichester, England, 1990).\\
$^{58}$D. Bao, S.-S. Chern and Z. Shen, \textit{An Introduction to
Riemann-Finsler Geometry} (Sprin\-ger-Verlag, New York, 2000).\\
$^{59}$Z. Shen, \textit{Differential Geometry of Spray and Finsler
Spaces.} (Kluwer, 2001).\\
$^{60}$R. Santamar\'{\i}a, ``Invariantes diferenciales de las
estructuras casi-bi\-paracomplejas y el problema de equivalencia,
Ph.D. thesis, Universidad de Cantabria (Santander, Spain), 2002.\\
$^{61}$B. Ammann, R. Lauter and V. Nistor, Int. J. Math. and Math.
Sciences (IJMMS) {\bf 4},  161 (2004).

\end{document}